\input amstex
\magnification=1200
\documentstyle{amsppt}
\NoRunningHeads
\NoBlackBoxes
\define\CA{\Cal A}

\define\CK{\Cal K}
\define\BC{\bold C}
\define\BO{\bold O}
\define\Fg{\frak g}
\define\Fv{\frak v}
\define\LC{\operatorname{LC}}
\define\GCR{\operatorname{\Gamma\Cal R}}
\define\Lie{\operatorname{\Cal L\Cal I\Cal E}}
\define\Abel{\operatorname{\Cal A\Cal B\Cal E\Cal L}}
\define\Add{\operatorname{\Cal A\Cal D\Cal D}}
\define\Mor{\operatorname{Mor}}
\define\End{\operatorname{End}}
\define\Ind{\operatorname{Ind}}
\define\HS{\operatorname{\Cal H\Cal S}}
\define\sltwo{\operatorname{\frak s\frak l}(2,\Bbb C)}
\topmatter
\title
Hidden symmetries and categorical representation theory.
\endtitle
\author Denis V. Juriev
\endauthor
\affil
``Thalassa Aitheria'' Research Center for Mathematical Physics and
Informatics,\linebreak
ul.Miklukho-Maklaya 20-180, Moscow, 117437, Russia.\linebreak
E-mail: denis\@juriev.msk.ru.\linebreak
\ \linebreak
\ \linebreak
\ \linebreak
November, 30, 1996\linebreak
q-alg/9612026\linebreak
\ \linebreak
\endaffil
\abstract The interrelations between the inverse problems of the
representation theory and the categorical representation theory are
discussed.
\endabstract
\endtopmatter
\document
\ \linebreak

This short article is devoted to the interrelations between the inverse
problems of the representation theory (see [1-4]) and the categorical
representation theory. The necessity of their analysis was clearly explicated
in the papers [3,4], where some non-standard representation theories appeared.
It seems that the machinery of the categorical representation theory is very
natural to work with some inverse problems of the representation theory.
Ideologically, it is truly remarkable that an analysis of concrete inverse
problems allows to enlight ``dark corners'' of the foundations of the theory.

The article is organized as a sequence of topics just as the previous
ones [2,3] and may be considered as their continuation devoted to
a specific subject. However, the concrete topics are prefaced by an
introduction, in which a general formalism of the categorical representation
theory is exposed. The paper contains three technically simple but
important theorems with straightforward and self--evident proofs, which are,
therefore, omitted.

\head Introduction: Elements of the categorical representation theory
\endhead

We shall consider the representations of classes of objects, which
constitute a category, which will be called {\it the ground category}.
The categorical aspects of the standard representation theory were discussed
in [5]. Some categorical generalizations were described in [6]. However, we
shall formulate the most abstract settings, which are necessary for our
purposes.

\definition{Definition 1A} {\it A representation theory\/} for the
ground category $\CA$ is a con\-tra\-va\-riant functor $R$ from the category
$\CA$ to the category $\Abel$ of all small abelian categories.
\enddefinition

Sometimes one should consider the category $\Add$ of all small additive
categories instead of $\Abel$. However, we shall consider the least
category for simplicity.

Often the ground category has some good properties, e.g. that for any
finite family of objects there exists their coproduct, which coincides
with their product. Such situation is realized for Lie algebras, Lie groups,
finite groups, associative algebras, Hopf algebras and many other structures.
However, the isotopic pairs (see f.e.[1:\S2.2]) and the most of other
algebraic pairs do not form a category of such type. For the ground
category $\CA$, in which products and coproducts of finite number of objects
exist and coincide, we shall claim in the definition of representation theories
that an associative family of imbeddings $R(a)\times R(b)\hookrightarrow
R(a+b)$ ($a$ and $b$ are any objects of the ground category $\CA$) is defined.
Such representation theories will be called {\it quasitensorial}.

\remark{Remark 1} If an object $a$ of the ground category $\CA$ admits a
coassociative mo\-no\-mor\-phism $\varepsilon$ into $a+a$ then $R(a)$ is
a tensor category iff the representation theory $R$ is quasitensorial.
\endremark

There exist non--quasitensorial representation theories even for the
well--known categories of the represented objects, e.g. general
$\HS$--projective representations of Lie algebras [3] or unitary
$\HS$--pseudorepresentations of Lie groups [4] are out of this class.

\definition{Definition 1B} A representation theory for the ground category
$\CA$ is called {\it ho\-mo\-mor\-phic\/} iff there exists a subcategory $\CA_0$ of
$\CA$ ({\it the target subcategory\/}) such that for any object $a$ of $\CA$
the category $R(a)$ may be identified with the category $\Mor(a,\CA_0)$ of all
(equivalence classes of) morphisms from $a$ to the objects of the category
$\CA_0$.
\enddefinition

For instance, theories of all linear, projective, unitary representations of
Lie groups are homomorphic. Note that the target category $\CA_0$ is always
an additive subcategory of the ground category $\CA$.

\definition{Definition 1C} A representation theory for the ground
category $\CA$ is called {\it hid\-den\-ly homomorphic\/} iff there exists a
homomorphic representation theory $R'$ for a category $\CK$ and a functor
(multi-valued as a rule) $\varrho:\CA\mapsto\CK$ such that $R=R'\circ\varrho$.
\enddefinition

Below we shall consider some examples and general constructions of the
hiddenly homomorphic representation theories for the ground category $\Lie$ of
the Lie al\-geb\-ras, which are not homomorphic, and describe their
interpretations in terms of the categorical representation theory, thus,
providing the least by the further elaboration of its details.

\head Topic One: The composite representation theories
\endhead

\definition{Definition 2 {\rm [3]}}

{\bf A.} A linear space $\Fv$ is called {\it a Lie composite\/}
iff there are fixed its subspaces $\Fv_1,\ldots \Fv_n$ ($\dim\Fv_i>1$)
supplied by the compatible structures of Lie algebras. Com\-pa\-ti\-bi\-li\-ty
means that the structures of the Lie algebras induced in $\Fv_i\cap\Fv_j$ from
$\Fv_i$ and
$\Fv_j$ are the same. The Lie composite is called {\it dense\/} iff
$\Fv_1\uplus\ldots\uplus\Fv_n=\Fv$ (here $\uplus$ means the sum of linear
spaces). The Lie composite is called {\it connected\/} iff for all $i$ and $j$
there exists a sequence $k_1,\ldots k_m$ ($k_1=i$, $k_m=j$) such that
$\Fv_{k_l}\cap\Fv_{k_{l+1}}\ne\varnothing$.

{\bf B.} {\it A representation\/} of the Lie composite $\Fv$ in the space $H$
is the linear mapping $T:\Fv\mapsto\End(H)$ such that $\left.T\right|_{\Fv_i}$
is a representation of the Lie algebra $\Fv_i$ for all $i$.

{\bf C.} Let $\Fg$ be a Lie algebra. A linear mapping $T:\Fg\mapsto\End(H)$
is called {\it the composed representation\/} of $\Fg$ in the linear space
$H$ iff there exists a set $\Fg_1,\ldots,\Fg_n$ of the Lie subalgebras of
$\Fg$, which form a dense connected composite and $T$ is its representation.
\enddefinition

Reducibility and irreducibility of representations of the Lie composites are
defined in the same manner as for Lie algebras. One may also formulate a
superanalog of the Definition 1. The set of representations of the fixed Lie
composite is closed under the tensor product and, therefore, may be supplied by
the structure of a tensor category. The theory of the composed
representations is evidently nonhomomorphic but hiddenly homomorphic theory.

The examples of the Lie composites and their representations were exposed
in [3]. The composed representations of the Witt algebra by the tensor
operators of spin 2 ($q_R$--conformal symmetries) in the Verma modules over
the Lie algebra $\sltwo$ were constructed in [3], too. They are representations
of the dense connected Witt composite and generate a tensor subcategory in the
category of all composed representations of the Witt algebra.

Below we shall formulate the general categorical setting for the construction
of the composed representations.

\definition{Definition 3A} Let $\CA$ be a topologized ground category
(i.e.supplied by a structure of the Grothendieck topology [7]). Let $R$
be a representation theory for $\CA$. {\it The composed representation
theory\/} $\BC(R)$ for $\CA$ may be constructed in the following manner.
Let $a$ be an object of the ground category $\CA$ and $S=(s_1,s_2,\ldots s_n)$
($s_i\in\Mor(a_i,a)$) be a cover of $a$ then the objects of the category
$\BC(R)(a)$ consists of all data $(b_1,b_2,\ldots b_n)$, $b_i\in R(a_i)$
such that for any object $c$ and monomorphisms $f\in\Mor(c,a)$ and
$f_i\in\Mor(c,a_i)$ ($f=s_i\circ f_i$) the equality ({\it the composite
glueing rule\/})
$$R(f_i)^*(b_i)=R(f_j)^*(b_j)$$
holds. The morphisms in $\BC(R)(a)$ are defined in the same manner.
\enddefinition

For any representation theory $R$ the composite representation theory
$\BC(R)$ is a sheaf of abelian categories over the topologized ground category
$\CA$ [7]. It is a sheaf canonically constructed from the pre--sheaf $R$
over the topologized ground category $\CA$ (note that the representation theory
for the topologized ground category $\CA$ is just a pre--sheaf over it).

\proclaim{Theorem 1} The composed representations of Lie algebras form a
composed representation theory $\BC(R)$, where $R$ is a standard
representation theory of Lie algebras (the covers of the Lie algebras are
defined by the dense connected Lie composites).
\endproclaim

Note that the Grothendieck topology of the Theorem 1 differs from the usual
one.

\remark{Remark 2} If $R$ is the standard representation theory then the
theory $\BC(R)$ is hiddenly homomorphic, the category $\CK$ is one of the
Lie composites, the category $\CK_0$ consists of Lie algebras $\End(H)$ for
all linear spaces $H$. i.e. just the same as for a homomorphic standard
represntation theory. However, if $R$ is a general representation theory
$\BC(R)$ is not obligatory hiddenly homomorphic.
\endremark

I suspect that the concept of the hidden homomorphicity of the composite
representation theories may be somehow understood in terms of the topos
theory [7].

\remark{Remark 3} $\BC(\BC(R))=\BC(R)$.
\endremark

\head Topic Two: The overlay representation theories
\endhead

The main disadvantage of the composed representation theory is clearly
explicated on the examples of the composed representations of the Witt algebras
by the hidden infinite dimensional ($q_R$--conformal) symmetries in the Verma
modules over the Lie algebra $\sltwo$. First, the tensor product of a finite
number of these irreducible composed representations is irreducible.
This fact contradicts to the na{\"\i}ve intuition. Second, the hidden symmetries
do not form any representation themselves whereas intuitively they should form
the adjoint representation.

So one needs some generalization of the composed representations.
Let us define the operator Lie composites $\LC(H)$ as the sets of subspaces
$\End(H_i)$ in the spaces $\End(H)$ ($H=H_1+\ldots+H_m$) with the natural
structures of Lie algebras.

\definition{Definition 4}

{\bf A.} {\it An overlay representation\/} of the Lie composite $\Fv$ in the
space $H$ is the homomorphism $T$ of $\Fv$ into the operator Lie composite
$\LC(H)$.

{\bf B.} Let $\Fg$ be a Lie algebra. A linear mapping $T:\Fg\mapsto\End(H)$
is called {\it the overlay composed representation\/} (or simply {\it overlay
representation\/}) of $\Fg$ in the linear space $H$ iff there exists a set
$\Fg_1,\ldots,\Fg_n$ of the Lie subalgebras of $\Fg$, which form a dense
connected composite and $T$ is its overlay representation.
\enddefinition

\remark{Remark 4} The overlay representations of any Lie algebra $\Fg$
form a tensor category.
\endremark

The overlay representations solve the previously described difficulties.

\proclaim{Theorem 2} The tensor operators of spin 2 in the Verma modules
$V_h$ over the Lie algebra $\sltwo$ form an overlay representation of the Witt
algebra, which are subrepresentations of $\End(V_h)$.
\endproclaim

\remark{Remark 5} The tensor operators of any natural spin $n$ in the Verma
modules $V_h$ over the Lie algebra $\sltwo$ form overlay representations of
the Witt algebra, which are subrepresentations of $\End(V_h)$.
\endremark

Let us formulate the natural categorical settings for the construction
of overlay representations.

\definition{Definition 3B} Let $\CA$ be a topologized ground category. Let $R$
be a homomorphic representation theory for $\CA$ with the target subcategory
$\CA_0$ supplied by the Grothendieck topology induced from $\CA$. {\it
The overlay representation theory\/} $\BO(R)$ for $\CA$ may be constructed
in the following manner. Let $a$ be an object of the ground category $\CA$ and
$S=(s_1,s_2,\ldots s_n)$ ($s_i\in\Mor(a_i,a)$) be a cover of $a$ then the
objects of the category $\BO(R)(a)$ consists of all data $(r_1,r_2,\ldots
r_n)$, $r_i\in\Mor(a,b_i)$, $b_i$ are objects of the target subcategory
$\CA_0$, which form a cover of the object $b$ of the same subcategory by the
monomorphisms $t_i\in\Mor(b_i,b)$, such that for any subobject $(c;p]$ of $a$
($p\in\Mor(c,a)$) the equality ({\it the overlay glueing rule\/})
$$r_i((a_i;s_i]\cap(c;p])\cap(b_j;t_j]=(b_i;t_i]\cap r_j((a_j;s_j]\cap(c;p])$$
holds. The morphisms in $\BO(R)(a)$ are defined in the same manner.
\enddefinition

However, I do not know a definition of the overlay representation theory
$\BO(R)$ for the representation theory $R$, which is not homomorphic.
Note that $\BO(R)$ is not a sheaf of abelian categories over $\CA$ in
general, and I do not know an abstract sheaf theoretical characterization
of the overlay representation theories.

\remark{Remark 6} The overlay representation theories $\BO(R)$ being defined
for the ho\-mo\-mor\-phic representations theories $R$ are hiddenly homomorphic.
\endremark

\proclaim{Theorem 3} The overlay representations of Lie algebras form
an overlay rep\-re\-sen\-tation theory $\BO(R)$, where $R$ is a standard
representation theory of Lie algebras (the covers of the Lie algebras are
defined by the dense connected Lie composites and the target subcategory
$\CA_0$ consists of all Lie algebras $\End(H)$).
\endproclaim

Note that the Grothendieck topology of the Theorem 3 differs from the usual
one.

\remark{Remark 7}
\roster
\item"--" If $R$ is a homomorphic representation theory for the ground category
$\CA$ then for any object $a$ of $\CA$ the category $\BC(R)(a)$ is a
subcategory of $\BO(R)(a)$.
\item"--" $\BC(\BO(R))=\BO(R)$.
\endroster
\endremark

\remark\nofrills{Problems:}
\roster
\item"--" To formulate (if possible) the definition of the overlay
representation theories $\BO(R)$ for the nonhomomorphic representation
theories $R$ (cf. a comment after the Definition 3B).
\item"--" To give a sheaf--theoretic description of the overlay representation
theories (cf. the same comment).
\item"--" To explicate the categorical or sheaf--theoretic relation between the
composite and overlay representation theories.
\item"--" To formulate the general categorical version of the induction
procedure (of the functors $\Ind_a^b:R(a)\mapsto R(b)$, $a$ is a subobject of
the object $b$ of the ground category) in the representation theory [5] and
adapt it to the composite and overlay representation theories.
\item"--" To adapt (if possible) the Tannaka--Krein theory [5] (see also [8])
to the overlay representations.
\endroster
\endremark

\remark{Remark 8} For some ground categories $\CA$ and some representation
theories $R$ (e.g. finite dimensional representations of finite groups and
reductive Lie algebras or unitary representations of compact groups) one may
construct the Grothendieck group (ring) $\Gamma(R(a))$ of virtual
representations ($a$ is any object of $\CA$) from the abelian (tensor)
category $R(a)$ [5] and, therefore, the pre--sheaves $\GCR$ of abelian groups
(rings) over $\CA$ (note that the Grothendieck rings $\Gamma(R(\cdot))$ are
supplied by numerous operations such as symmetric and exterior degrees, Adams
operations and arbitrary ``polynomial operations'' related to the irreducible
representations of the symmetric groups $S_n$ [5]). The transition to the
composite representation theory supply the topologized ground category $\CA$
by the sheaf $\BC(\GCR)$ over it, so one may consider the cohomologies
$H^*(\CA,\BC(\GCR))$ of $\CA$ with coefficients in $\BC(\GCR)$. The additional
operations supply the cohomologies $H^*(\CA,\BC(\GCR))$ by a sophisticated
algebraic structure.
\endremark

\head Conclusions
\endhead

Thus, a categorical framework for the non--standard representation theories,
which appear in the analysis of hidden symmetries [3,4], is briefly described.
The relations of the abstract categorical constructions to the concrete
inverse problems of the representation theory are explicated. Open questions
are formulated. It seems that as the abstract categorical formalism as its
appearance in the context of the analysis of hidden symmetries may be useful
for the understanding of some general representation theoretic aspects of the
control theory [9] and applications of the representation theory to the
concrete problems of (classical and quantum) controlled systems as well as for
many other subjects.

\Refs
\roster
\item"[1]" Juriev D., An excursus into the inverse problem of representation
theory [in Russian]. Report RCMPI-95/04 (August 1995) [e-version:
mp\_arc/96-477].
\item"[2]" Juriev D., Topics in hidden symmetries. I-IV. E-prints:
hep-th/9405050, q-alg/9610026, q-alg/9611003, q-alg/9611019.
\item" [3]" Juriev D., Topics in hidden symmetries. V. E-print:
funct-an/9611003.
\item"[4]" Juriev D., On the infinite dimensional hidden symmetries. I.
Infinite dimensional geometry of $q_R$-conformal symmetries. E-print:
funct-an/9612004.
\item"[5]" Kirillov A.A., Elements of the representation theory. Springer,
1976; An introduction into the representation theory and harmonic analysis
[in Russian]. Current Probl.Math. Fundam.Directions, V.22, VINITI, Moscow,
1988.
\item"[6]" Geronimus A.Yu., Grothendieck topology and representation theory
[in Russian]. Funkts. anal.i ego prilozh. 5(3) (1971) 22-31.
\item"[7]" Th\'eorie des topos et cohomologie \'etale des sch\'emas. Berlin,
Heidelberg, New York, 1972-1973; Cohomologie \'etale. Berlin, Heidelberg,
New York, 1977.
\item"[8]" Naimark M.A., Normed rings [in Russian]. Moscow, Nauka, 1968.
\item"[9]" Juriev D., Topics in hidden algebraic structures and
infinite--dimensional dynamical symmetries of controlled systems, in
preparation.
\endroster
\endRefs
\enddocument